\documentclass[a4paper,11pt,UKenglish]{article}

\usepackage{amsmath,amssymb, xcolor, listings, amsthm}
\usepackage{mathtools}
\usepackage{enumitem}
\usepackage{url}
\usepackage{hyperref}
\usepackage{fullpage}

\usepackage[linesnumbered,boxed,vlined]{algorithm2e}	
\SetAlCapSkip{0.5em}

\usepackage{graphicx}
\usepackage{pstricks}
\usepackage{caption}
\usepackage{pst-plot, pst-node, pst-text, pst-tree}
\usepackage{rotating}
\usepackage{epic, eepic}
\usepackage{color}
\usepackage{indentfirst}
\usepackage{mathrsfs}
\usepackage{manfnt}

\hypersetup{colorlinks=true, linkcolor=black, citecolor=blue, urlcolor=blue}
\theoremstyle{definition}

\theoremstyle{plain}
\newtheorem{theorem}{Theorem}[section]
\newtheorem{prop}[theorem]{Proposition}
\newtheorem{lemma}[theorem]{Lemma}
\newtheorem{cor}[theorem]{Corollary}
\DeclarePairedDelimiter{\abs}{\lvert}{\rvert} 
\newcommand{\multinom}[2]{\left(\!\!\binom{#1}{#2}\!\!\right)}		
\newcommand{\inlinemultinom}[2]{\left(\!\binom{#1}{#2}\!\right)}	
\newcommand{\cupdot}{\mathbin{\mathaccent\cdot\cup}}    

\author{Lapo Cioni\thanks{Dipartimento di Matematica e Informatica ``U. Dini'',
Universit\`a degli Studi di Firenze, Firenze, Italy. {
\tt\ \{lapo.cioni,luca.ferrari\}@unifi.it}. Member of the INdAM research group GNCS; partially supported by the 2020 INdAM-GNCS project "Combinatoria delle permutazioni, delle parole e dei grafi:
algoritmi e applicazioni".}\and Luca Ferrari$^*$}

\title{Preimages under the Queuesort algorithm}

\begin{document}

\maketitle

\begin{abstract}
Following the footprints of what have been done with the algorithm \texttt{Stacksort}, 
we investigate the preimages of the map associated with a slightly less well known algorithm, called \texttt{Queuesort}.
After having described an equivalent version of \texttt{Queuesort}, we provide a recursive description of the set of all preimages of a given permutation,
which can be also translated into a recursive procedure to effectively find such preimages. 
We then deal with some enumerative issues. 
More specifically, we investigate the cardinality of the set of preimages of a given permutation, showing that all cardinalities are possible, except for 3.
We also give exact enumeration results for the number of permutations having 0,1 and 2 preimages.
Finally, we consider the special case of those permutations $\pi$ whose set of left-to-right maxima is the disjoint union of a prefix and a suffix of $\pi$:
we determine a closed formula for the number of preimages of such permutations, which involves two different incarnations of ballot numbers, 
and we show that our formula can be expressed as a linear combination of Catalan numbers.   
\end{abstract}

\section{Introduction}

\texttt{Stacksort} is a classical and well-studied algorithm that attempts to sort an input permutation by (suitably) using a stack.
It has been introduced and first investigated by Knuth \cite{K} and West \cite{W},
and it is one of the main responsible for the great success of the notion of \emph{pattern} for permutations.
Among the many research topics connected with \texttt{Stacksort},
a very interesting one concerns the characterization and enumeration of preimages of the associated map, which is usually denoted with $s$
(so that $s(\pi )$ is the permutation which is obtained after performing \texttt{Stacksort} on $\pi$).
More specifically, given a permutation $\pi$, what is $s^{-1}(\pi )$? How many permutations does it contain?
These questions have been investigated first by Bousquet-Melou \cite{BM}, and more recently by Defant \cite{D1} and Defant, Engen and Miller \cite{DEM}.

\bigskip

In the present paper we address the same kind of problems for a similar sorting algorithm.
Suppose to replace the stack with a queue in \texttt{Stacksort}. What is obtained is a not so useful algorithm,
whose associated map is the identity (and so, in particular, the only permutations that it sorts are the identity permutations).
However, if we allow one more operation, namely the bypass of the queue, the resulting algorithm (which we call \texttt{Queuesort}) is much more interesting.
This is not a new algorithm, and some properties of it can be found scattered in the literature \cite{B,JLV,M,T}. However, to the best of our knowledge, the problem of studying preimages under the map associated with \texttt{Queuesort}  (similarly to what have been done for \texttt{Stacksort}) has never been considered. Our aim is thus to begin the investigation of this kind of matters, with a particular emphasis on enumeration questions.

\bigskip

Our paper is organized as follows. In Section \ref{prel} we provide basic notions and terminology concerning \texttt{Queuesort}, as well as some fundamental preliminary results, that will be heavily used throughout the paper.
In particular, we describe an alternative version of \texttt{Queuesort}, that is an algorithm which is different from \texttt{Queuesort} but completely equivalent to it. 
Our algorithm acts directly on the input permutation by moving some of its elements,
without any reference to the queue. As a matter of fact, throughout all the paper we will work with this alternative version of \texttt{Queuesort} in order to obtain all our results.
Section \ref{rec_char} contains a characterization of preimages of a given permutation $\pi$, from which a recursive algorithm to generate all preimages of $\pi$ is obtained.
One of the main tool to get such a recursive description is a special decomposition of permutations, which we have called \emph{LTR-max decomposition}: it will prove to be extremely useful also in subsequent sections, and a key ingedient for most of our achievements.
Results concerning enumeration can be found in Section \ref{enumeration}, where we start by showing how the number of preimages of a permutations only depends on the position of its LTR maxima (rather than their values). We then provide some specific results concerning the exact number of permutations having few preimages, and we show that all cardinalities are possible for the set of preimages of a given permutation, except for 3. Finally, we analyze those permutations $\pi$ in which the set of all LTR maxima is the disjoint union of a prefix and a suffix of $\pi$; for them, we are able to determine a closed formula for the number of preimages (which is in general rather complicated), involving well-known quantities such as ballot numbers, multinomial coefficients and Catalan numbers.
Section \ref{conclusion} contains some hints for further work.

\section{Preliminary notions and results}\label{prel}

Given a permutation $\pi=\pi_1 \pi_2 \cdots \pi_n$, the algorithm \texttt{Queuesort}
 attempts to sort $\pi$ by using a queue in the following way: scan $\pi$ from left to right and, called $\pi_i$ the current element,
\begin{itemize}
	\item if the queue is empty or $\pi_i$ is larger than the last element of the queue, then $\pi_i$ is inserted to the back of the queue;
	\item otherwise, compare $\pi_i$ with the first element of the queue, then output the smaller one.
\end{itemize}
When all the elements of $\pi$ have been processed, pour the content of the queue into the output. For a more formal description, see Algorithm \ref{queuesort}.

\begin{algorithm}
	$Queue:=\emptyset$\;
	$i:=1$\;
	\While{$i\leq n$}
	{
		\If{$Queue=\emptyset$ or $BACK(Queue) < \pi_i$}
		{
			$\textnormal{execute}$ Q\;
			$i:=i+1$\;
		}
		\Else
		{
			\While{$FRONT(Queue) < \pi_i$}{$\textnormal{execute}$ O\;}
			$\textnormal{execute}$ B\;
			$i:=i+1$\;
		}
	
	}
	\While{$Queue \neq \emptyset$}
	{$\textnormal{execute}$ O\;}
	\caption{\texttt{Queuesort} ($Queue$ is the queue; $FRONT(Queue)$ is the current first element of the queue; $BACK(Queue)$ is the current last element of the queue; operations Q,B,O are \textquotedblleft insert into the queue", \textquotedblleft bypass the queue", \textquotedblleft place the first element of the queue into the output", respectively;  $\pi=\pi_1 \cdots \pi_n$ is the input permutation).}\label{queuesort}
\end{algorithm}

There is a natural function $q$ associated with \texttt{Queuesort}, which maps an input permutation $\pi$ into the permutation $q(\pi )$ that is obtained by performing \texttt{Queuesort} on $\pi$.

The set of permutations that can be sorted using \texttt{Queuesort} (simply called \emph{sortable} from now on) can be suitably characterized by means of the notion of pattern. Recall that, given two permutations $\sigma =\sigma_1 \sigma_2 \cdots \sigma_k ,\tau =\tau_1 \tau_2 \cdots \tau_n$, with $k\leq n$, we say that $\sigma$ is a \emph{pattern} of $\tau$ whenever there exist indices $1\leq i_1 <i_2 <\cdots <i_k\leq n$ such that $\sigma$ is order-isomorphic to $\tau_{i_1}\tau_{i_2}\cdots \tau_{i_k}$ (i.e., for all $l,m$, $\sigma_l < \sigma_m$ if and only if $\tau_{i_l} <\tau_{i_m}$). The set of all permutations avoiding a pattern $\sigma$ is denoted $Av(\sigma )$, whereas $Av_n (\sigma )$ is the subset of $Av(\sigma )$ consisting of permutations of length $n$.
As it was shown by Tarjan \cite{T}, a permutation is sortable if and only if it avoids the pattern 321. In terms of its geometric structure, a 321-avoiding permutation is a permutation that can be expressed as the shuffle of two increasing subpermutations. In other words, if $id_n =12\cdots n$ is the identity permutation of length $n$, $q^{-1}(id_n )=Av_n (321)$. Moreover, as it is well-known, $|q^{-1}(id_n )|=|Av_n (321)|=C_n$, where $C_n =\frac{1}{n+1}\binom{2n}{n}$ is the $n$-th Catalan number.

\bigskip

In order to study preimages of the map $q$, the first tool we develop is an effective description of the behavior of \texttt{Queuesort} directly on the input permutation. This will prove very useful in the sequel. Before starting, we need one more definition. An element $\pi_i$ of a permutation $\pi$ is called a \emph{left-to-right maximum} (briefly, \emph{LTR maximum}) when it is larger than every element to its left (that is, $\pi_i >\pi_j$, for all $j<i$).

Let $\pi$ be a permutation of length $n$, and denote with $m_1 ,m_2 ,\ldots ,m_k$ its LTR maxima, listed from left to right. Thus, in particular, $m_k =n$. Then $q(\pi )$ is obtained from $\pi$ by moving its LTR maxima to the right according to the following instructions:
\begin{itemize}
	\item for $i$ running from $k$ down to 1, repeatedly swap $m_i$ with the element on its right, until such an element is larger than $m_i$.
\end{itemize}

For instance, if $\pi =21543$, then there are two LTR maxima, namely 2 and 5; according to the above instructions, $\pi$ is thus modified along the following steps: $\mathbf{2}1\mathbf{5}43\rightsquigarrow \mathbf{2}1435\rightsquigarrow 12435$, and so $q(21543)=12435$.

\begin{center}
\begin{minipage}{.20\textwidth}
\centering
\includegraphics[trim=7.5cm 20cm 7.5cm 4cm, width=\textwidth]{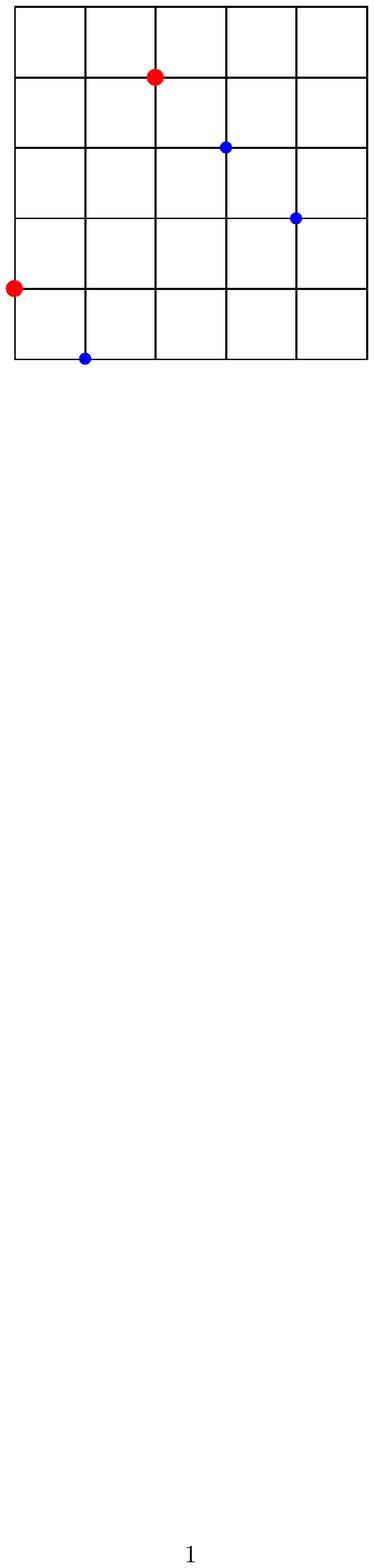}
\end{minipage}
$\overrightarrow{\quad}$
\begin{minipage}{.20\textwidth}
\centering
\includegraphics[trim=7.5cm 20cm 7.5cm 4cm, width=\textwidth]{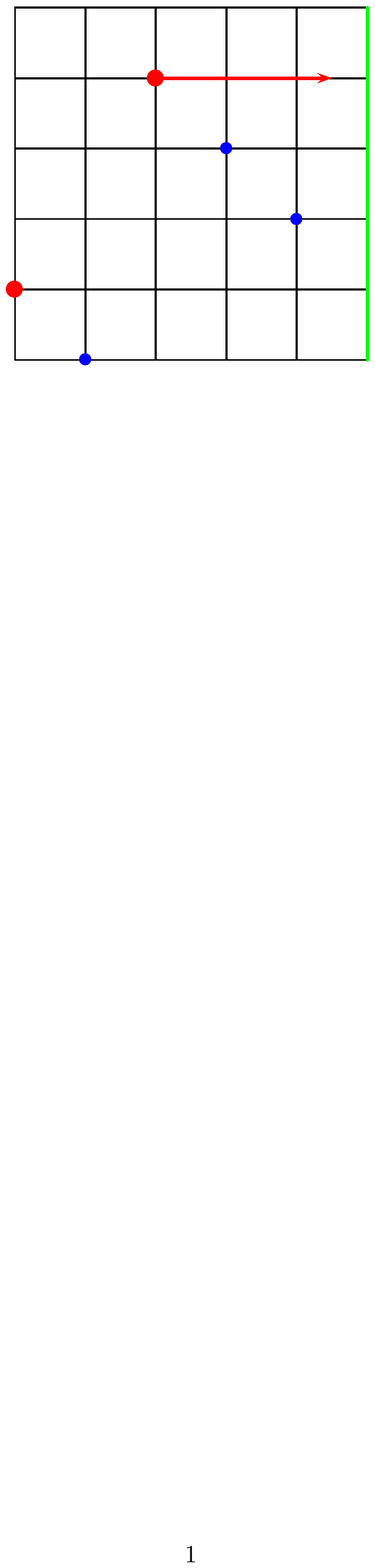}
\end{minipage}
$\overrightarrow{\quad}$
\begin{minipage}{.20\textwidth}
\centering
\includegraphics[trim=7.5cm 20cm 7.5cm 4cm, width=\textwidth]{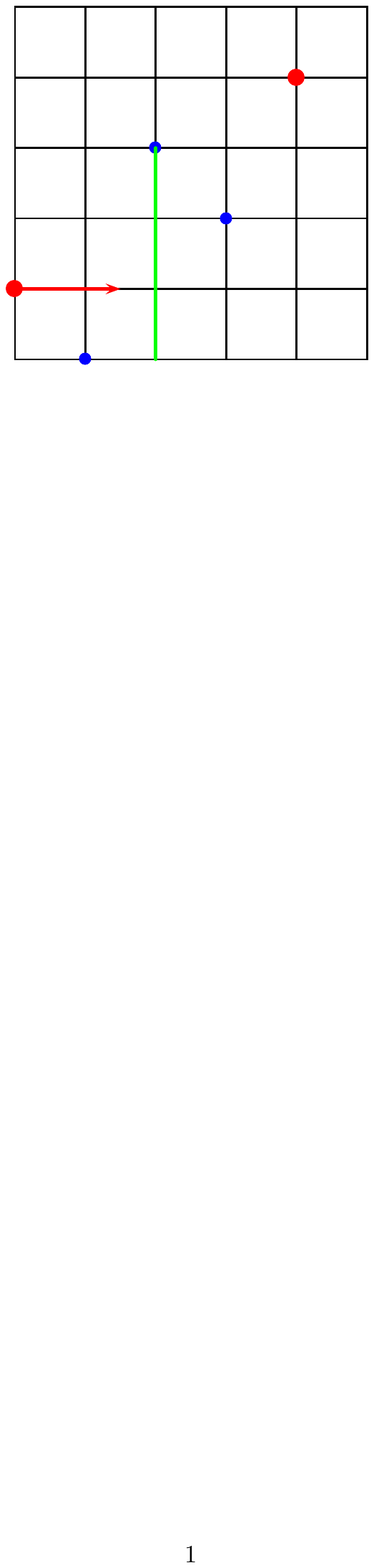}
\end{minipage}
$\overrightarrow{\quad}$
\begin{minipage}{.20\textwidth}
\centering
\includegraphics[trim=7.5cm 20cm 7.5cm 4cm, width=\textwidth]{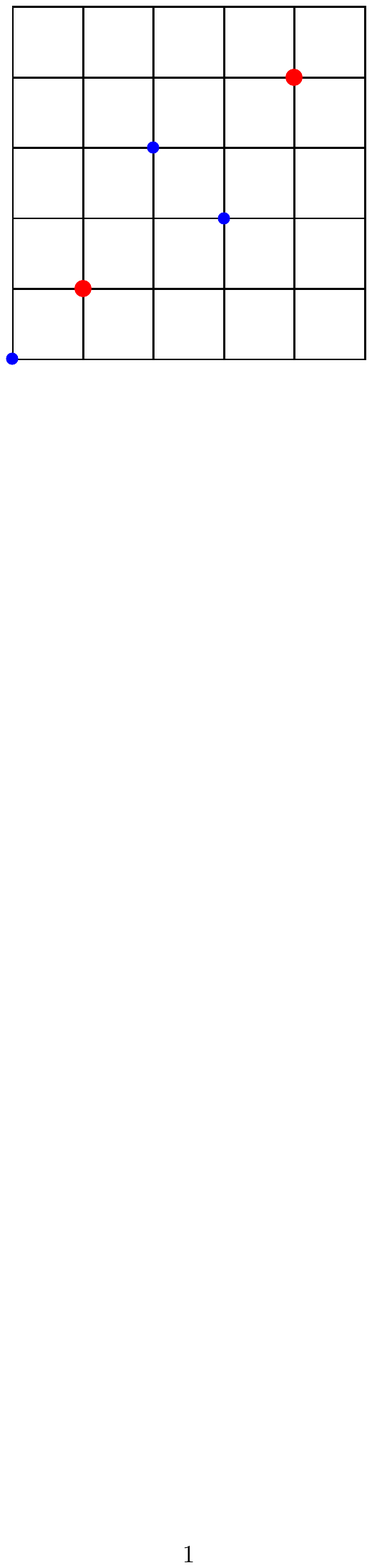}
\end{minipage}
\end{center}

The above alternative description of \texttt{Queuesort} is, as a matter of fact, a different algorithm which is however equivalent to \texttt{Queuesort}
(meaning that, starting from a given input, it returns exactly the same output). The proof of this fact is easy, and it relies on the fact that the elements of the input permutations that enters the queue are precisely its LTR maxima.
In the rest of the paper, with some abuse of terminology, whenever we will declare that we perform the algorithm \texttt{Queuesort},
we will in fact perform our equivalent algorithm.

\bigskip

An immediate consequence of this alternative description of \texttt{Queuesort} is our first result, characterizing permutations having nonempty preimage.

\begin{prop}\label{zero}
Given $\pi =\pi_1 \pi_2 \cdots \pi_n$, we have that $q^{-1}(\pi )\neq \emptyset$ if and only if $\pi_n =n$.
\end{prop}

\begin{proof}
Suppose first that $q^{-1}(\pi )\neq \emptyset$, and let $\sigma \in q^{-1}(\pi )$. From the above description of \texttt{Queuesort}, it follows immediately that the last element of $q(\sigma)=\pi$ is $n$. Conversely, suppose that $\pi_n =n$, and define $\sigma =n\pi_1 \cdots \pi_{n-1}$. It is easy to check that $q(\sigma )=\pi$, and so $\sigma \in q^{-1}(\pi )\neq \emptyset$.
\end{proof}

In closing this section, we collect a few notations that will be used throughout the whole paper.

We denote with $S_n$ the set of all permutations of length $n$. Given $\pi \in S_n$ written in one-line notation, we will denote with $\pi_k$ the $k$-th element of $\pi$, so that $\pi =\pi_1 \pi_2 \cdots \pi_k$. The identity permutation of length $n$ will be denoted $id_n$. To conclude, we observe that, in the sequel, it will often happen that, given a permutation $\pi$, we construct new permutations by removing one or more elements from $\pi$. Strictly speaking, the resulting sequence is not a permutation, but it can be considered as such after suitably rescaling its elements. In what follows we will always assume that such a rescaling is performed, so we will omit to explicitly recall it in every single case. 

\section{A recursive characterization of preimages}\label{rec_char}

The key ingredient to state our recursive characterization of preimages is a suitable decomposition of permutations,
which is based on the notion of LTR maximum.
Given a permutation $\pi$, we decompose it as $\pi=M_1 P_1 M_2 P_2 \cdots M_{k-1} P_{k-1} M_k$,
where the $M_i$'s are all the maximal sequences of contiguous LTR maxima of $\pi$ (and the $P_i$'s collect all the remaining elements). This decomposition will be called the \emph{LTR-max decomposition} of $\pi$, and it will prove extremely useful in the sequel.
In particular, all the $P_i$'s are nonempty, and $M_i$ is nonempty for all $i\neq k$.
Moreover, $m_i =|M_i |$ denotes the length of $M_i$, and analogously $p_i =|P_i |$ denotes the length of $P_i$, for all $i$.
Sometimes we will also need to refer to the last element of $M_i$, which will be denoted $\mu_i$.
Also, in some cases we will use $N$ and $R$ in place of $M$ and $P$, respectively. In order to avoid repeating the same things several times, the above notations for the LTR-max decomposition of a permutation will remain fixed throughout all the paper.

In order to determine the preimages of a given permutation $\pi$,
we can exploit the description of \texttt{Queuesort} illustrated in the previous section,
observing that every preimage of $\pi$ can be obtained by moving suitable elements of $\pi$ to the left, until they reach a suitable position.
The next proposition (whose easy proof is left to the reader) is recorded for further reference, and says that the only elements of $\pi$ that can be moved are LTR maxima of $\pi$.

\begin{prop}\label{LTR}
Suppose that $\pi =q(\sigma )$. Then each LTR maximum of $\sigma$ is a LTR maximum of $\pi$ as well.
\end{prop}

The next result is instead less trivial, and gives a necessary condition for a LTR maximum of $\pi$ to be movable in order to get preimages.

\begin{prop}\label{not-LTR}
Suppose that $\pi =q(\sigma )$, with $\pi ,\sigma \in S_n$. Let $\pi_i$ be a LTR maximum of $\pi$, for some $i<n$.
If $\pi_{i+1}$ is not a LTR maximum of $\pi$, then $\pi_i$ is not a LTR maximum of $\sigma$.
\end{prop}

\begin{proof}
Suppose, by contradiction, that $\pi_i$ is a LTR maximum of $\sigma$.
According to our description of \texttt{Queuesort}, consider the instant when $\pi_i$ is moved to the right.
When $\pi_i$ reaches its final position, it is clear that all the elements $\pi_{i+1},\pi_{i+2},\ldots ,\pi_n$ are already in their final positions in $\pi$
(since, when a LTR maximum is moved to its final position, no further elements can overtake it).
This means that, when $\pi_i$ reaches its final position, it is immediately followed by $\pi_{i+1}$.
However, by hypothesis, $\pi_{i+1}$ is not a LTR maximum of $\pi$, so necessarily $\pi_i >\pi_{i+1}$, hence $\pi_i$ should overtake $\pi_{i+1}$,
which gives a contradiction.
\end{proof}

As a consequence of the above proposition, whenever $\pi_i$ is not followed by a LTR maximum in $\pi$,
we cannot move it to the left in order to get a preimage of $\pi$.

\bigskip

We are now ready to start our description of all the preimages of a given permutation $\pi$.
Our approach consists of splitting the set of preimages into two disjoint subsets, each of which can be described in a recursive fashion.
For the rest of the section, we will assume that $\pi=M_1 P_1 M_2 P_2 \cdots M_{k-1}P_{k-1}M_k \in S_n$ is a permutation having at least one preimage.
We know that this is the same as saying that the last element of $\pi$ is $n$, or equivalently that $M_k$ is nonempty.
Moreover, we will suppose that $\pi$ is not the identity permutation (whose preimages are well understood).

\begin{prop}\label{first}
Suppose that $|M_k |\geq 2$. Denote with $\pi'$ the permutation obtained from $\pi$ by removing $n$, and let $\sigma'$ be any preimage of $\pi'$.
Let $\sigma'=N_1 R_1 \cdots N_{s-1}R_{s-1}N_s$ be the LTR-max decomposition of $\sigma'$.
Then every permutation obtained from $\sigma'$ by inserting $n$ in any position to the right of $N_{s-1}$ is a preimage of $\pi$.
\end{prop}

\begin{proof}
Let $\sigma$ be obtained from $\sigma'$ by inserting $n$ somewhere to the right of $N_{s-1}$. Our goal is to show that $q(\sigma )=\pi$.
Observe that, when we perform \texttt{Queuesort} on $\sigma'$, the elements of $N_s$ are not moved to the right.
Therefore, when we perform \texttt{Queuesort} on $\sigma$, first of all $n$ is moved to the rightmost position,
then (since $n$ was to the right of $N_{s-1}$) the algorithm performs the same operations it would perform on $\sigma'$.
As a consequence, $q(\sigma )=q(\sigma')n=\pi' n=\pi$, as desired.
\end{proof}

\begin{prop}\label{second}
Let $M_{k-1}'$ (respectively, $M_k '$) be the (possibly empty) sequence obtained by removing the last element $\mu_{k-1}$ (respectively, $n$)
from $M_{k-1}$ (respectively $M_k$).
Denote with $\sigma$ any preimage of the permutation $\rho =M_1 P_1 \cdots M_{k-2}P_{k-2}M_{k-1}'n$.
Then the permutation $\tau$ obtained by concatenating $\sigma$ with $\mu_{k-1}P_{k-1}M_k '$ is a preimage of $\pi$.
\end{prop}

\begin{proof}	
When performing \texttt{Queuesort} on $\tau$, first of all $n$ is moved to the rightmost position.
After that, observe that the element $\mu_{k-1}$ is larger than all the elements on its left.
Hence, all the elements of $\sigma$ will not overtake it during the execution of the algorithm.
Therefore, at the end, the output permutation is obtained by concatenating $q(\sigma )$, $\mu_{k-1}P_{k-1}M_k '$ and $n$,
i.e. $q(\tau )=\pi$, as desired.
\end{proof}

The previous two propositions told us how to find some preimages of $\pi$.
The next two propositions aim at showing that these are the only preimages of $\pi$.

\begin{prop}\label{third}
Suppose that $|M_k |\geq 2$. Denote with $\pi'$ the permutation obtained from $\pi$ by removing $n$.
Let $\sigma =N_1 R_1 \cdots N_{s-1}R_{s-1}N_s$ be any preimage of $\pi$ such that
there exists an element of $M_k$ different from $n$ which belongs neither to $R_{s-1}$ nor to $N_s$.
Denote with $\sigma'$ the permutation obtained from $\sigma$ by removing $n$.
Then $\sigma'$ is a preimage of $\pi'$.
\end{prop}

\begin{proof}
Looking at the position of $n$ in $\sigma$, we can distinguish two cases.

If $n\in N_s$, then $\sigma =\sigma' n$, and clearly $q(\sigma' n)=q(\sigma )=\pi =\pi' n$, hence $q(\sigma')=\pi'$, as desired.

On the other hand, suppose that $n\notin N_s$. This means of course that $N_s =\emptyset$, and $n$ is the rightmost element of $N_{s-1}$.
Let $D=\{ \gamma \in M_k \, |\, \gamma \neq n, \gamma \notin N_s , \gamma \notin R_{s-1}\}$, and set $\alpha =\max D$.
Notice that $\alpha$ indeed exists, since our hypothesis implies that $D\neq \emptyset$.

When we perform \texttt{Queuesort} on $\sigma$, first of all the element $n$ is moved to the rightmost position.
After that, we now want to show that the next element to be moved is $\alpha$.

First, we notice that $\alpha$ is in fact a LTR maximum of $\sigma$ (hence it is moved by \texttt{Queuesort}).
Indeed, if this were not the case, $\alpha$ would be to the left of $N_{s-1}$ in $\sigma$,
and so it could never overtake the elements of $N_{s-1}$ and $R_s$, against the fact that $\alpha \in M_k$.
Our next aim is to prove that there cannot exist elements between $\alpha$ and $n$ in $\sigma$ that can be moved by \texttt{Queuesort}.
By contradiction, suppose that $\beta$ is such an element.
Then necessarily $\beta$ is a LTR maximum of $\sigma$, and it belongs neither to $R_{s-1}$ (since $R_{s-1}$ does not contain LTR maxima)
nor to $N_s$ (because elements in $N_s$ are already at the end of the permutation and so they are not moved by \texttt{Queuesort}).
Moreover, $\beta \neq n$ (by our assumption) and $\beta >\alpha$ (since $\beta$ is a LTR maximum to the right of $\alpha$ in $\sigma$).
Finally, $\beta \in M_k$, since $\alpha \in M_k$ is a LTR maximum of $\pi$,
hence $\beta$ ($>\alpha$) cannot be in a position to the left of $M_k$ in $\pi$.
We can thus conclude that $\beta \in D$, but this gives a contradiction, since $\max D=\alpha <\beta \in D$.

To conclude the proof, we now want to show that $q(\sigma' )n=q(\sigma)$ (since obviously $q(\sigma )=\pi =\pi' n$, hence $q(\sigma' ) =\pi'$). A moment's thought should convince that, to this aim, it is enough to prove that there exist no element after $n$ in $\sigma$ which is larger than $\alpha$ and is the first element of a descent. For a contradiction, suppose that $\omega$ is such an element. Then $\alpha$ cannot overtake $\omega$ (when \texttt{Queuesort} is performed on $\sigma$), and this implies that $\alpha \notin M_k$ (since there is at least one element following $\alpha$ in $\pi$ which is smaller than $\alpha$), which is false.
\end{proof}

\begin{prop}\label{fourth}
Let $M_{k-1}'$ (respectively, $M_k '$) be the (possibly empty) sequence obtained by removing the last element $\mu_{k-1}$ (respectively, $n$)
from $M_{k-1}$ (respectively $M_k$).
Denote with $\sigma =N_1 R_1 \cdots N_{s-1}R_{s-1}N_s$ a preimage of $\pi$ such that
all the elements of $M_k$ other than $n$ belong to $R_{s-1}N_s$.
Then $\sigma$ is the concatenation of
a preimage of $M_1 P_1 \cdots M_{k-2}P_{k-2}M_{k-1}'n$ and $\mu_{k-1}P_{k-1}M_k '$.
\end{prop}

\begin{proof}
We start by showing that $n$ is the last element of $N_{s-1}$ (and so $N_s =\emptyset$).
Since $P_{k-1}\neq \emptyset$, the element immediately after $\mu_{k-1}$ in $\pi$ is not a LTR maximum of $\pi$.
Thus, by Proposition \ref{not-LTR}, $\mu_{k-1}$ is not a LTR maximum of $\sigma$.
Denote with $\nu$ any LTR-maximum of $\sigma$ preceding $\mu_{k-1}$ such that $\nu >\mu_{k-1}$.
As a consequence of Proposition \ref{LTR}, all the LTR maxima of $\sigma$ larger than $\mu_{k-1}$ must belong to $M_k$.
However, since $\nu$ precedes $\mu_{k-1}$ (and $\mu_{k-1}$ is not a LTR maximum of $\sigma$), $\nu$ cannot belong to $R_{s-1}N_s$.
Therefore we can conclude that $\nu =n$,
hence $n$ is the last element of $N_{s-1}$ and $\sigma =N_1 R_1 \cdots N_{s-1}R_{s-1}$.

Now decompose $\sigma$ as $\sigma =\rho \tau$, where the first element of $\tau$ is $\mu_{k-1}$.
Performing \texttt{Queuesort} on $\sigma$, we have first that $n$ is moved to the last position.
After that, we observe that all the LTR maxima of $\sigma$ that are moved are smaller than $\mu_{k-1}$.
In fact, all LTR maxima of $\sigma$ larger than $\mu_{k-1}$ are also LTR maxima of $\pi$ (by Proposition \ref{LTR}),
therefore they are to the right of $\mu_{k-1}$ in $\sigma$, i.e. they belong to $\tau$.
Thus we have that $M_1 P_1 M_2 P_2 \cdots M_{k-1}P_{k-1}M_k =\pi =q(\sigma )=\theta \tau n$,
for a suitable $\theta$, hence $\tau =\mu_{k-1}P_{k-1}M_k '$.

Finally,
since all the LTR maxima of $\sigma$ (other than $n$) which are moved by \texttt{Queuesort} are smaller than $\mu_{k-1}$,
they cannot overtake $\mu_{k-1}$, and so $\pi =q(\sigma )=q(\rho )'\tau n$
(where, as usual, $q(\rho )'$ denotes the permutation obtained from $q(\rho )$ by removing $n$).
Thus $q(\rho )=q(\rho )'n=M_1 P_1 \cdots M_{k-2}P_{k-2}M_{k-1}'n$, which concludes the proof.
\end{proof}

The last four propositions allow us to state the announced recursive description of all preimages of a given permutation.

\begin{theorem}\label{recurrence}
Let $\pi =M_1 P_1 \cdots M_{k-1}P_{k-1}M_k \in S_n$, with $M_k \neq \emptyset$, and suppose that $\pi$ is different from the identity permutation.
A permutation $\sigma \in S_n$ is a preimage of $\pi$ if and only if exactly one of the following holds:
\begin{itemize}
\item $\sigma =\tau \mu_{k-1}P_{k-1}M_k '$, where, with a little abuse of notation, $\tau \in q^{-1}(M_1 P_1 \cdots M_{k-2}P_{k-2}M_{k-1}'n)$ denotes a preimage of $M_1 P_1 \cdots M_{k-2}P_{k-2}M_{k-1}'n$ and $M_{k-1}'$, $M_k '$ are defined as in Proposition \ref{fourth};
\item if $\pi'$ is defined (as in Proposition \ref{third}) by removing $n$ from $\pi$
and $\sigma' =N_1 R_1 \cdots N_{s-1}R_{s-1}N_s$ is a preimage of $\pi'$,
then $\sigma$ is obtained by inserting $n$ in one of the positions to the right of $N_{s-1}$.
\end{itemize}
\end{theorem}

\begin{proof}
Notice that, for a given $\pi$, exactly one of the above two cases hold:
indeed, in the second case there is at least one element of $M_k$ other than $n$ which is moved, whereas this does not happen in the first case.
In particular, in the second case a preimage exists only if $|M_k |\geq 2$ (otherwise no preimage of $\pi'$ can exist).
Thus the theorem is a consequence of Propositions \ref{first} to \ref{fourth}.
\end{proof}

Thanks to the above theorem, we can also describe a recursive algorithm to determine all preimages of a given permutation
$\pi =M_1 P_1 \cdots M_{k-1}P_{k-1}M_k$:

\begin{itemize}

\item if $\pi$ is the identity permutation, then the preimages of $\pi$ are precisely the 321-avoiding permutations (of the same length);

\item otherwise

\begin{itemize}

\item compute all the preimages of $M_1 P_1 \cdots M_{k-2}P_{k-2}M_{k-1}'n$ and concatenate them with $\mu_{k-1}P_{k-1}M_k '$;

\item if $|M_k |\geq 2$, then compute all the preimages $N_1 R_1 \cdots N_{s-1}R_{s-1}N_s$ of $\pi'$
and insert $n$ in each of the positions to the right of $N_{s-1}$.

\end{itemize}

\end{itemize}

\emph{Example.}\quad To illustrate the above theorem, we now compute all the preimages of the permutation 23145. The various steps of the procedure are depicted in the figure below:

\begin{center}
\includegraphics[trim=5.5cm 18cm 5.5cm 5cm,width=0.5\textwidth]{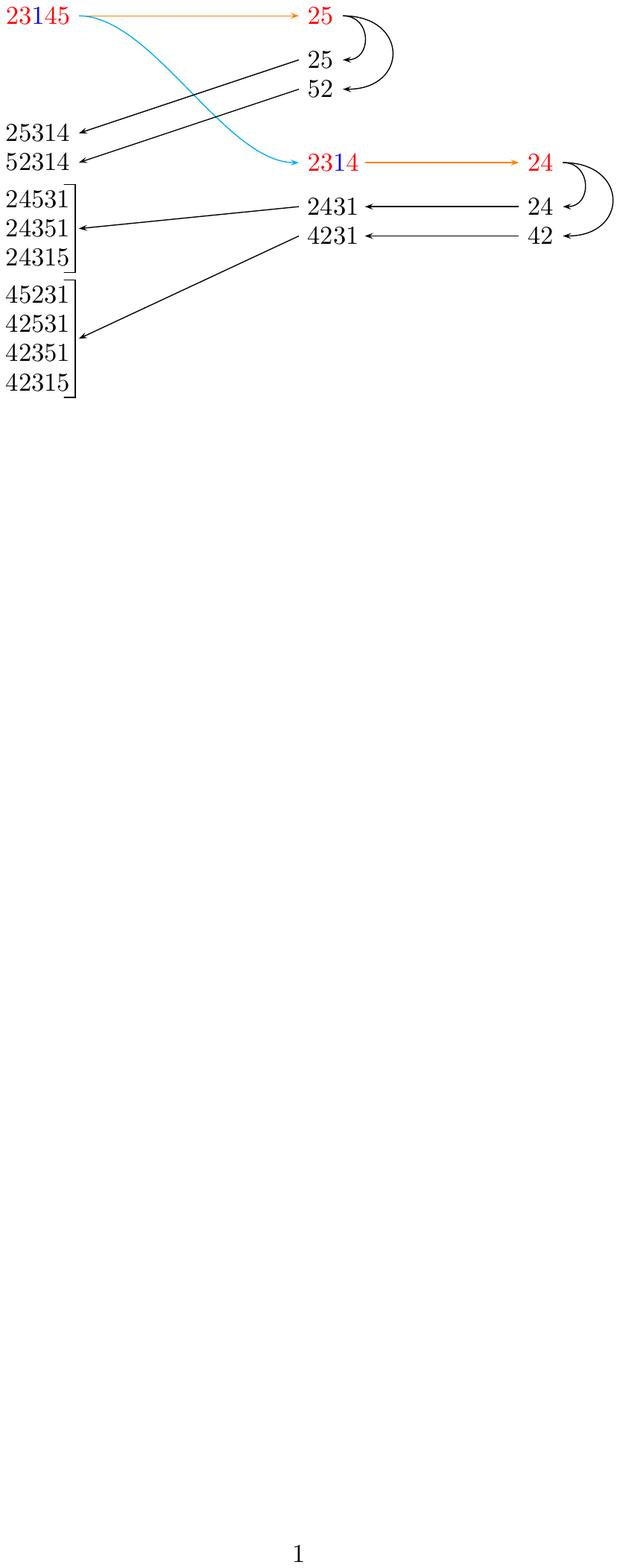}
\end{center}

\section{Enumerative results}\label{enumeration}

The main goal of the present section is to count how many preimages under \texttt{Queuesort} a given permutation has.
We begin by proving some preliminary facts, which are however relevant enough to deserve a proper presentation.

\bigskip

Our first achievement is an important feature of the \texttt{Queuesort} algorithm:
the number of preimages of a permutation depends only on the positions of its LTR maxima, and not on their values.
The proof of this property will stem from a more general fact, that reveals an interesting feature of \texttt{Queuesort}.

Given $\pi \in S_n$, denote with $LTR(\pi )$ the set of the positions of the LTR maxima of $\pi$,
i.e. $LTR(\pi )=\{ i\leq n\, |\, \pi_i \textnormal{ is a LTR maximum of } \pi \}$.

\begin{lemma}\label{goingback}
Let $\pi, \sigma \in S_n$, with $\sigma$ obtained from $\pi$ by moving some of its LTR maxima to the left, from positions $i_1,\dots , i_h$ to positions $i'_1,\dots , i'_h$, where $i_1<\cdots <i_h$, $i'_1<\cdots <i'_h$, $i'_j<i_j$ for every $j=1,\dots ,h$. Then, for every $j=1,\ldots ,h$, the position of $\sigma_{i'_j}$ in $q(\sigma )$ is at least $i_j$.
\end{lemma}

\begin{proof}
We start by observing that all the elements of $\pi$ that are moved to the left must remain LTR maxima also in $\sigma$.
We now perform \texttt{Queuesort} on $\sigma$ and describe what happens to the LTR maxima in positions $i'_1,\ldots ,i'_h$. The element $\sigma_{i'_h}$ is the first one which is considered. Since all the elements between positions $i'_h$ and $i_h$ are smaller that $\sigma_{i'_h}$ (by construction of $\sigma$), its final position in $q(\sigma )$ will be at least $i_h$. Moving on, suppose that all the elements $\sigma_{i'_{j+1}},\ldots ,\sigma_{i'_h}$ have already been processed, and their final positions are at least $i_{j+1},\ldots ,i_h$, respectively. When \texttt{Queuesort} considers $\sigma_{i'_j}$, the only elements between positions $i'_j$ and $i_j$ that could possibly be greater than $\sigma_{i'_j}$ are the previously moved LTR maxima, whose positions are however strictly greater than $i_j$ (since $i_j <i_{j+1}$). Therefore, also the final position of $\sigma_{i'_j}$ in $q(\sigma )$ is at least $i_j$. Repeating the argument for all the moved LTR maxima gives the thesis.     
\end{proof}

\begin{lemma}
\label{oltrepassare}
Let $\pi=\pi_1\cdots\pi_n = M_1 P_1 \cdots M_{k-1} P_{k-1} M_k \in S_n$ and let $\sigma$ be a preimage of $\pi$. As usual, let $\mu_i$ be the last element of $M_i$. Then, for all $i=1,\ldots ,k-1$, there exists a LTR maximum of $\pi$ which is to the left of $\mu_i$ in $\sigma$ and to the right of $\mu_i$ in $\pi$.	
\end{lemma}

\begin{proof}
By Proposition \ref{not-LTR}, $\mu_i$ is not a LTR maximum of $\sigma$. Thus there must be a LTR maximum of $\sigma$ to the left of $\mu_i$ which is greater than it. Such an element cannot be to the left of $\mu_i$ in $\pi$, since $\mu_i$ is a LTR maximum of $\pi$, and this gives the thesis.     
\end{proof}

\begin{cor}
\label{lastblock}
Using the notations of the above lemma, the largest element of $\sigma$ moved by \texttt{Queuesort} belongs to $M_k$.
\end{cor}

\begin{theorem}
\label{LTRsubset}
Let $\rho ,\pi \in S_n$. If $LTR(\rho )\subseteq LTR(\pi )$, then $|q^{-1}(\rho )|\leq |q^{-1}(\pi )|$.  	
\end{theorem}

\begin{proof}
Define a map $f:q^{-1}(\rho )\rightarrow q^{-1}(\pi )$ as follows. Given $\tau \in q^{-1}(\rho )$, suppose that, performing \texttt{Queuesort} on $\tau$, the LTR maxima of $\tau$ in positions $i'_1 <\cdots <i'_h$ are moved to the right to positions $i_1 <\cdots <i_h$, respectively. Define $f(\tau )=\sigma$ as the permutation obtained from $\pi$ by moving the elements in positions $i_1 ,\ldots ,i_h$ to the left to positions $i'_1 ,\ldots ,i'_h$, respectively.

We now show that $f$ is well defined, i.e. that indeed $\sigma \in q^{-1}(\pi )$. Since $LTR(\rho )\subseteq LTR(\pi )$, $\pi_{i_1} ,\ldots ,\pi_{i_h}$ are LTR maxima of $\pi$, so, by Lemma \ref{goingback}, performing \texttt{Queuesort} on $\sigma$ moves them at least to their original positions (possibly to the right of them).
First of all, we prove that no further elements of $\sigma$ are moved by \texttt{Queuesort}. 
By contradiction, let $\sigma_{a'}$ be a LTR maximum of $\sigma$ that is moved by \texttt{Queuesort}, with $a'\neq i_j$ for every $j$. Let $a$ be such that $\pi_a=\sigma_{a'}$. Since, by construction of $\sigma$, the set of elements preceding $\pi_a$ in $\pi$ is a subset of the elements preceding $\sigma_{a'}$ in $\sigma$, we have that $\pi_a$ is a LTR maximum of $\pi$. If we now look at $\rho$, we have two possibilities concerning the index $a$: either $a\notin LTR(\rho)$ or $a\in LTR(\rho)$.

If $a\notin LTR(\rho)$, set $b=max \{ c<a \mid \rho_c \text{ is a LTR maximum of } \rho \}$. By Lemma \ref{oltrepassare} applied to $\rho$, there exists a $\bar{t}$ such that $i'_{\bar{t}}<b<i_{\bar{t}}$. In particular, by definition of $b$ we have $i'_{\bar{t}}<b<a<i_{\bar{t}}$, and so $\pi_a<\pi_{i_{\bar{t}}}$ (because $i_{\bar{t}}\in LTR(\rho) \subseteq LTR(\pi)$). Since in $\sigma$ the element $\pi_a$ is to the right of the element $\pi_{i_{\bar{t}}}$, that is in position $i'_{\bar{t}}$ (as a consequence of the way $\sigma$ is obtained from $\pi$), we have that $\sigma_{a'}=\pi_a$ is not a LTR maximum of $\sigma$, which is a contradiction.\\
If instead $a\in LTR(\rho)$, we can assume that there exists no $\bar{t}$ such that $i'_{\bar{t}}<a<i_{\bar{t}}$, because otherwise the same argument of the previous case would lead to the contradiction that $\sigma_{a'}$ is not a LTR maximum of $\sigma$. We now consider $\rho$. In particular, denoting with $M_j$ the block of LTR maxima containing $\rho_a$, let $b$ be the index of the rightmost element of $M_j$. Applying Lemma \ref{oltrepassare}, we find that there exists $\bar{u}$ such that $i'_{\bar{u}}<b<i_{\bar{u}}$. By the above assumption, we have $a<i'_{\bar{u}}$. This implies that, when we apply \texttt{Queuesort} to $\sigma$, the element in position $i'_{\bar{u}}$ will move at least to position $i_{\bar{u}}$, and the element $\sigma_{a'}=\pi_a$ does not change position because all the elements between $\pi_a$ and $\pi_b$ (included) do not move. This is however in contrast with the hypothesis that $\sigma_{a'}$ is moved by \texttt{Queuesort}.

We have thus shown that the only elements in $\sigma$ moved by \texttt{Queuesort} are those in positions $i'_1,\ldots, i'_h$. It remains to prove that they are moved to their original positions $i_1,\ldots,i_h$.
We first observe that $\sigma_{i'_h}$ goes back to position $i_h$ because, applying Corollary \ref{lastblock} to $\rho$, we obtain that (when applying \texttt{Queuesort}) $\tau_{i'_h}$ goes back to the last block of LTR maxima of $\rho$, and so $\sigma_{i'_h}$ goes back to the last block of LTR maxima of $\pi$ as well, and this means that it cannot move further.
Now suppose by contradiction that there exists some $j<h$ such that the element $\sigma_{i'_j}$ is moved by \texttt{Queuesort} to a position strictly to the right of $i_j$. Let $\bar{j}$ be the maximum of such $j$'s. Since $\rho=q(\tau)$, $\rho$ has a LTR maximum in position $i_{\bar{j}}+1$, so the same is true for $\pi$. Now observe that, in $\sigma$, \texttt{Queuesort} moves $\sigma_{i'_{\bar{j}}}$ to a position strictly to the right of $i_{\bar{j}}$, so at that point of the execution of \texttt{Queuesort} the element $\pi_{i_{\bar{j}}+1}$ is not in position $i_{\bar{j}}+1$. Thus it must be to the right of that position (it cannot be to the left by Lemma \ref{goingback}), which is against the maximality of $\bar{j}$.

We have thus shown that the map $f$ is well defined. To conclude, we observe that it is also injective, because each preimage of $\rho$ is uniquely determined by the $i_j$'s and $i'_j$'s.
\end{proof}

\begin{cor}
If two permutations $\pi$ and $\rho$ have their LTR maxima in the same positions, then they have the same number of preimages.
\end{cor}
\begin{proof}
By hypothesis we have $LTR(\pi)=LTR(\rho)$, so the previous theorem implies that $\abs{q^{-1}(\pi)}=\abs{q^{-1}(\rho)}$.
\end{proof}

Another relevant consequence of Theorem \ref{LTRsubset} is the following proposition, which reveals to be a useful tool in several circumstances.

\begin{prop}
\label{noLTRisol}
Let $\pi =M_1 P_1 \cdots M_{k-1}P_{k-1}M_k \in S_n$, with $|M_i|=|\{ \mu_i \} |=1$ for a given $i\neq 1,n$.
Let $\rho =N_1 R_1 \cdots N_{i-1} R'_i N_{i+1} \cdots N_{k-1}R_{k-1}N_k \in S_n$, with $|M_j |=|N_j |$, $|P_j |=|R_j |$, for all $j\neq i$, $|P_{i-1}M_i P_i |=|R'_i |$ and such that $R'_i$ does not contain any LTR maximum of $\rho$.
Then $|q^{-1}(\rho )|=|q^{-1}(\pi )|$.
\end{prop}

\begin{proof}
We start by observing that $LTR(\rho)\subset LTR(\pi)$, so, by Theorem \ref{LTRsubset}, we have $|q^{-1}(\rho )|\leq |q^{-1}(\pi )|$.

We now show that $|q^{-1}(\rho )|\geq |q^{-1}(\pi )|$. In fact, due to Proposition \ref{not-LTR}, $\mu_i$ cannot be moved when looking for a preimage of $\pi$. So, given a preimage of $\pi$, we can construct a preimage of $\rho$ as described in the proof of Theorem \ref{LTRsubset} and show that it is indeed a preimage of $\rho$ using similar arguments.
Therefore $\rho$ has at least as many preimages as $\pi$, thus giving the desired inequality.
\end{proof}

In other words, the previous proposition says that, in some sense, the presence of isolated LTR maxima does not affect the number of preimages.

\bigskip

We now provide some results concerning permutations with a given number of preimages.
We already know (Proposition \ref{zero}) that $\pi\in S_n$ has no preimages if and only if its last element is different from $n$. Therefore, setting $Q_n ^{(k)}=\{ \pi \in S_n \, |\, |q^{-1}(\pi )|=k\}$ and $q_n ^{(k)}=|Q_n ^{(k)}|$, we have that $Q_n ^{(0)}=\{ \pi \in S_n  \, |\, \pi_n \neq 0\}$ and $q_n ^{(0)}=(n-1)!\cdot (n-1)$. The next propositions deal with $Q_n ^{(1)}$ and $Q_n ^{(2)}$.

\begin{prop}
For all $n$, we have $Q_n ^{(1)}=\{  \pi \in S_n \, |\, \pi_n =n$ and $\pi$ does not have two adjacent LTR maxima$\}$. As a consequence, 
\[
q_n ^{(1)}=(n-1)!\cdot \sum_{i=0}^{n-1}\frac{(-1)^i}{i!},
\] 
that is the $(n-1)$-th derangement number (sequence A000166 in \cite{Sl}).	
\end{prop}

\begin{proof}
Suppose first that $\pi \in S_n$ does not have any two adjacent LTR maxima and its last element is $n$.
Then certainly $|q^{-1}(\pi )|\neq 0$.
Using repeatedly Proposition \ref{noLTRisol}, we can assert that $|q^{-1}(\pi )|=|q^{-1}(\rho )|$,
where $\rho$ is any permutation of length $n$ whose only LTR maxima are in positions 1 and $n$.
Thus, in particular, $|q^{-1}(\pi )|=|q^{-1}((n-1)(n-2)\cdots 21n)|$,
and it is clear (by a direct computation, or by invoking Theorem \ref{recurrence}) that the last quantity is equal to 1.

On the other hand, suppose that the last element of $\pi$ is $n$ (otherwise, of course, $\pi$ has no preimages),
but there exists $i$ such that $\pi_i ,\pi_{i+1}$ are both LTR maxima.
We could now invoke Theorem \ref{recurrence} to find two distinct preimages of $\pi$.
However we prefer to explicitly describe two such preimages, since it is so simple.
One preimage of $\pi$ can be obtained by moving $n$ at the beginning of the permutation;
in other words, $\pi_n \pi_1 \pi_2 \cdots \pi_{n-1}$ is a preimage of $\pi$ (this is trivial to verify).
Another preimage of $\pi$ can be obtained by placing $n$ between $\pi_i$ and $\pi_{i+1}$
and replacing $\pi_1 \cdots \pi_i$ with any of its preimages $\tau$
(notice that the set of preimages of $\pi_1 \cdots \pi_i$ is indeed not empty, since $\pi_i$ is a LTR maximum);
in other words, $\tau n\pi_{i+1}\cdots \pi_{n-1}$ is a preimage of $\pi$ (this is also quite easy to realize, and so left to the reader). The two above described preimages are indeed distinct, since the former starts with $n$ whereas the latter does not.

To conclude the proof, we recall the so called Foata's fundamental bijection \cite{FS},
which maps a permutation $\sigma$ written in one-line notation to the permutation in cycle notation
obtained by inserting a left parenthesis in $\sigma$ preceding every LTR maximum,
then a right parenthesis where appropriate.
Applying such a map to $Q_n (1)$ returns the set of permutations of length $n$ whose only fixed point is $n$,
which is clearly equinumerous with the set of derangements of length $n-1$, as desired.
\end{proof}

\begin{prop}
For all $n$, we have $Q_n ^{(2)}=\{  \pi \in S_n \, |\, \pi_n =n$ and $\pi$ does not have two adjacent LTR maxima except for the first two elements$\}$. As a consequence, $q_n ^{(2)}$ satisfies the recurrence relation
\begin{align*}
&q_{n+1} ^{(2)}=(n-1)q_{n}^{(2)}+(n-1)q_{n-1}^{(2)},\qquad n\geq 3,\\
&q_0 ^{(2)}=q_1 ^{(2)}=q_3 ^{(2)}=0,\quad q_2 ^{(2)}=1.
\end{align*}
\end{prop}

\begin{proof}
By an argument completely analogous to that of the previous proposition, we can show that any permutation ending with $n$ and whose only adjacent LTR maxima are the first two elements has exactly two preimages, that are obtained by moving $n$ into the first and into the second position, respectively.

On the other hand, if $\pi \in S_n$ does not have adjacent LTR maxima, we know by the previous proposition that $|q^{-1}(\pi )|=1$. Moreover, if $\pi$ contains two adjacent LTR maxima $\pi_i =\alpha $ and $\pi_{i+1}$, with $i\geq 2$, then $\pi$ has at least three preimages, obtained as follows: either move $n$ into the first position, or move both $\alpha$ and $n$ into the first two positions, or replace $\alpha$ with $n$ and move $\alpha$ into the first position.


We will now prove the recurrence relation. We can immediately see that the initial conditions hold.
First observe that, if we take a permutation in $Q_{n+1}^{(2)}$, then we can remove $n+1$, which is the last element, to obtain a permutation of length $n$ without any two adjacent LTR maxima except for the first two elements, and such that its last element is not $n$. We can now use Foata's fundamental correspondence \cite{FS} to get a bijection with the set $A_n$ of permutations of length $n$ with only one fixed point, which is smaller than the maximum of each of the remaining cycles. Thus, set $a_n =\abs{A_n}$, the recurrence for $q_{n}^{(2)}$ is equivalent to the following, which we will now prove:
\[
a_{n}=(n-1)a_{n-1}+(n-1)a_{n-2},\qquad n\geq 4.
\]

Given $\rho\in A_n$, its maximum $n$ necessarily belongs to a cycle with at least two elements. If such a cycle has exactly two elements, say $n$ and $\alpha$, with $1\le\alpha\le n-1$, then we can remove them to obtain a permutation of length $n-2$ belonging to $A_{n-2}$, thus obtaining the second summand of the r.h.s of the recurrence equation, i.e. $(n-1)a_{n-2}$.
If instead the cycle contains at least three elements, we have two distinct cases. If it contains at least one other element larger than the fixed point, then we can remove $n$ to obtain a permutation in $A_{n-1}$. Observe that, to invert this construction, given a permutation in $A_{n-1}$, there are $n-2$ possible positions to insert $n$ in order to get a permutation of $A_n$ (since $n$ can be placed before every element except for the fixed point).
Finally, if all the elements of the cycle containing $n$ are smaller than the fixed point, except for $n$, then we can proceed as follows.
Let $\alpha$ be the element following $n$ in its cycle (i.e. $\alpha =\rho_n$ in the one-line notation of $\rho$) and let $\beta$ be the fixed point. Remove $n$ from $\rho$, and switch $\alpha$ with $\beta$. By construction, we obtain a permutation of length $n-1$  with exactly one fixed point, which is smaller than the maximum of each cycle, and from which we can revert back to the starting permutation uniquely. The last two cases together account for the first summand of the r.h.s of the recurrence relation, i.e. $(n-1)a_{n-1}$.
\end{proof}

Sequence $q_n ^{(2)}$ starts $0,0,1,0,2,6,32,190\ldots$ and is essentially A055596 in \cite{Sl}. We thus deduce, for $n\geq 2$, the closed formula $q_n ^{(2)}=(n-1)!-2q_{n}^{(1)}$ (recall that $q_n ^{(1)}$ equals the $(n-1)$-th derangement number), as well as the exponential generating function
\[
\sum_{n\geq 0}q_n ^{(2)}\frac{x^n}{n!}=\frac{x(2-x-2e^{-x})}{1-x}.
\]

We have thus seen that there exist permutations having 0, 1 or 2 preimages. We now show (Propositions \ref{not3} and \ref{=3}) that there exist permutations having any number of preimages, except for 3. Before that, we need a preliminary result, which is of interest in its own.

\begin{lemma}\label{catalan}
Let $\pi =M_1 P_1 M_2$, with $|M_2|=1$. Then $|q^{-1}(\pi )|=C_{m_1}$, the $m_1$-th Catalan number.
\end{lemma}

\begin{proof}
By hypothesis $\pi =M'_1 \mu_1 P_1 n$, where $\mu_1$ is the last element of $M_1$.
Applying Theorem \ref{recurrence}, we obtain that all the preimages of $\pi$ are of the form $\tau \mu_1 P_1$, with $\tau \in q^{-1}(M'_1 n)$.
However, the permutation $M'_1 n$ is increasing, hence $|q^{-1}(M'_1 n)|=C_{m_1}$, which gives the thesis.  
\end{proof}

\begin{prop}\label{not3}
Given $n\geq 2$, let $\pi =n(n-1)(n-2)\cdots 21(n+2)(n+3)(n+1)(n+4)\in S_{n+4}$.
Then $|q^{-1}(\pi )|=n+2$.
\end{prop}

\begin{proof}
We will repeatedly use Theorem \ref{recurrence} to compute preimages.
First of all, observe that only the first case of such a theorem applies to $\pi$,
since the second-to-last element of $\pi$ is not $n+3$.
Hence, the preimages of $\pi$ are precisely those permutations of the form $\sigma (n+3)(n+1)$,
where $\sigma$ is any preimage of $n\cdots 1(n+2)(n+4)$.
Therefore our problem reduces to the enumeration of such $\sigma$.
Preimages of $n\cdots 1(n+2)(n+4)$ can again be computed using Theorem \ref{recurrence}.
In particular, the first case of that theorem gives rise to a single preimage, which is $(n+4)n\cdots 1(n+2)$.
Looking at the second case, we first have to compute the set of preimages of $n\cdots 1(n+2)$,
which is easily seen to consist of the single permutation $\tau =(n+2)n\cdots 1$ of length $n+1$
(notice that this holds since we are supposing that $n\geq 2$, so that $\tau$ is not an increasing permutation);
then we have to suitably insert $n+4$ into $\tau$, and this can be done in $n+1$ different ways,
namely by inserting $n+4$ to the right of each element of $\tau$.
Thus we have found that the total number of preimages of $\pi$ is precisely $n+2$, as desired.
\end{proof}

The previous proposition can be easily extended to the case $n=0$,
since it is immediate to check that $|q^{-1}(2314)|=2$.
However, it does not hold when $n=1$, as $|q^{-1}(13425)|=5$.
The next proposition shows that this is no accident.

\begin{prop}\label{=3}
There exists no permutation $\pi$ such that $|q^{-1}(\pi )|=3$.
\end{prop}

\begin{proof}
As usual, suppose that $\pi$ is decomposed as $\pi =M_1 P_1 \cdots M_{k-1}P_{k-1}M_k$.
If $k=1$, then $|q^{-1}(\pi )|$ is a Catalan number, which is certainly different from 3.
If $k\geq 2$, we distinguish two cases, depending on the cardinality of $M_2$.

\begin{itemize}

\item If $|M_2 |>1$, denote with $\alpha$ and $\beta$ the two largest elements of $M_2$, with $\alpha <\beta$, so that $M_2 =M''_2 \alpha \beta$.
Then the following four permutations are all preimages of $\pi$, and are all distinct 
(this is a consequence of Theorem \ref{recurrence}, or it can be checked by applying to each of them our equivalent version of \texttt{Queuesort}):
\begin{itemize}
\item $\beta M_1 P_1 M''_2 \alpha P_2 \cdots M_k$;
\item $\alpha \beta M_1 P_1 M''_2 P_2 \cdots M_k$;
\item $\alpha M_1 \beta P_1 M''_2 P_2 \cdots M_k$;
\item $\alpha M_1 \beta P_1 M''_2 \beta P_2 \cdots M_k$.
\end{itemize}

\item  If $|M_2 |=1$, then Proposition \ref{noLTRisol} tells us that the number of preimages of $\pi$ is the same as 
the number of preimages of any permutation $\sigma= N_1 R_1 \cdots N_{k-2}R_{k-2}N_{k-1}$ such that 
$|N_1 |=|M_1 |$, $|R_1 |=|P_1 M_2 P_2 |$ and $|N_i |=|M_{i+1} |$, $|R_i |=|P_{i+1} |$, for all $i\geq 2$.
We can iterate this argument until 
either the second maximal sequence of consecutive LTR maxima of the resulting permutation has cardinality $\geq 2$ 
or we obtain a permutation $\sigma$ of the form $\sigma =N_1 R_1 N_2$, with $|N_2 |=1$.
In the former case (which occurs when there exists $i\geq 2$ such that $|M_i |\geq 2$), we are in the situation described in the previous item point.
In the latter case (which occurs when $|M_i |=1$, for all $i\geq 2$), we can apply Lemma \ref{catalan},
thus obtaining again that $|q^{-1}(\pi )|$ is a Catalan number.
\end{itemize}
\vspace{-25pt}
\end{proof}

The last part of our paper is devoted to find an expression for the number of preimages of a generic permutation $\pi$ of the form $\pi =M_1 P_1 M_2$.
We have already proved (Lemma \ref{catalan}) that, when $|M_2 |=1$, we get a Catalan number.
In order to find a general formula, we need some preliminary work.

For any $n\ge 1$ and $2\le i\le n$, define $G_{n,i}=\{\pi\in S_n(321) \mid \pi_i \text{ is not a LTR maximum and, for every } j<i,\ \pi_j \text{ is a LTR maximum} \}$, and $g_{n,i}=\abs{G_{n,i}}$. Also define $G_{n,n+1}=\{ id_n \}$, $g_{n,n+1}=1$.

\begin{prop}
For every $n\ge 2$, the following recurrence holds:
\begin{equation}
g_{n,i}=
\begin{cases}
C_{n-1}, \quad &\text{if $i=2$,}\\
n-1, \quad &\text{if $i=n$,}\\
g_{n-1,i-1}+g_{n,i+1} \quad &\text{if $3\le i\le n-1$.}
\end{cases}
\end{equation}
\end{prop}
\begin{proof}
We observe that $G_{n,2}=\{\pi\in S_n(321) \mid \pi_2 = 1\}$, which is in bijection with $S_{n-1}(321)$ (by the removal of $\pi_2$). As a consequence, $g_{n,2}=C_{n-1}$.\\
Moreover, given $\pi \in G_{n,n}$, the first $n-1$ elements are in increasing order, and the last element can be chosen arbitrarily, except that it cannot be $n$. Thus $g_{n,n}=\abs{G_{n,n}}=n-1$.\\
To prove the recurrence relation, we describe a bijection $f$ from $G_{n,i}$ to $G_{n-1,i-1}\cupdot G_{n,i+1}$ (where $\cupdot$ denotes disjoint union). Given $\pi\in G_{n,i}$, we define $f(\pi)$ as follows:
\[
f(\pi)=
\begin{cases}
\pi_2\cdots\pi_n &\in G_{n-1,i-1},\quad \text{if $\pi_1=1$},\\
\pi_1\cdots\pi_{i-1}\pi_{i+1}\pi_i\pi_{i+2}\cdots\pi_n &\in G_{n,i+1},\quad \quad \text{if $\pi_1\neq1$ and $\pi_{i+1}>\pi_{i-1}$},\\
\pi_i\pi_1\cdots\pi_{i-1}\pi_{i+1}\cdots\pi_n &\in G_{n,i+1},\quad \quad \text{if $\pi_1\neq1$ and $\pi_{i+1}<\pi_{i-1}$}.
\end{cases}
\]

It is easy to check that $f$ is indeed a bijection by explicitly constructing its inverse.
\end{proof}

The above recurrence relation also extends to the case $i=n$ (and also $i=n+1$, by setting all undefined values to be $0$).

\begin{cor}\label{gi}
For every $n\ge 1$, $2\le i\le n+1$,
$g_{n,i}=\binom{2n-i+1}{n}\frac{i-1}{2n-i+1}$.
\end{cor}
\begin{proof}
Up to suitably rescaling the indices, sequence $g_{n,i}$ satisfies the same initial conditions and recurrence relation as sequence $A033184$ in \cite{Sl}. More specifically, $g_{n,i}=A033184(n,i-1)$. This gives the desired closed formula.
\end{proof}

If we represent sequence $A033184$ as a triangle, we obtain one of the so-called Catalan triangles. Its entries are sometimes called ballot numbers \cite{A}.

Before stating and proving the main result of this section, we need to introduce one more statistic on $321$-avoiding permutations. Namely, we define $b_{n,i} = \abs{\{\pi\in S_n(321)\mid\pi_i=n\}}$. Using the bijection that associates every permutation with its group-theoretic inverse, we can see that $b_{n,i} = \abs{\{\pi\in S_n(321)\mid\pi_n=i\}}$.

Next proposition shows that the $b_{n,i}$'s are another version of the ballot numbers.
\begin{prop}\label{bi}
The numbers $b_{n,i}$ correspond to sequence $A009766$ of \cite{Sl}, if we properly shift them. Namely, the element $b_{n,i}$ corresponds to the element of indices $n-1,i-1$ of sequence $A009766$.
\end{prop}

\begin{proof}
We will prove the statement by showing that the $b_{n,i}$'s satisfy the same recurrence relation (and initial conditions) as sequence A009766, namely $b_{n+1,i}=\sum_{j=1}^{i}b_{n,j}$ (and $b_{1,1}=1$, which is plainly true).
Let $\pi \in S_{n+1}(321)$ such that $\pi_{n+1}=i$. 
By removing $i$ (and rescaling), we get a permutation $\pi' \in S_n (321)$.
If $\pi'_n =j\neq n$, then necessarily $j<i$ (since $\pi$ avoids 321), and each such $\pi'$ gives uniquely a permutation in $S_{n+1}(321)$ ending with $i$ (by appending $i$ to $\pi'$).
On the other hand, if $j=n$, then $\pi'=M_1 P_1 \cdots M_{k-1}P_{k-1}M_k$ is such that $M_k$ is not empty, and the last element of $P_{k-1}$ is smaller than $i$, because otherwise $\pi$ would contain the pattern $321$. Thus, denoting with $\alpha$ the number of such permutations, we have
\[
b_{n+1,i}=\sum_{j=1}^{i-1}b_{n,j}+\alpha .
\]

To determine $\alpha$ we observe that, removing $M_k$ from $\pi'$, if $\pi'$ is not the identity permutation, we obtain a bijection with the set of permutations avoiding $321$ of any length $m<n$ whose last element $j$ is smaller than $i$ and different from $m$.
Hence, using induction (on $n$) and repeatedly applying the recurrence relation for the $b_{n,i}$'s recalled at the beginning of this proof, we get that
\[
\alpha=1+\sum_{m=2}^{n-1} \sum_{j=1}^{\min{(m-1,i-1)}} b_{m,j}=b_{n,i}.
\]
Indeed, the inner sum for a fixed $m$ is the same as the l.h.s. of the recurrence relation, except for the element $b_{m,m}$, which is obtained by the inner sum for $m-1$ (for $m=2$, we observe that $b_{2,2}=1$).
Summing up, we get:
\[
b_{n+1,i}=\sum_{j=1}^{i-1}b_{n,j}+b_{n,i}=\sum_{j=1}^{i}b_{n,j},
\]
as desired.
\end{proof}

\emph{Remark.}\quad Corollary \ref{gi} and Proposition \ref{bi} imply that $g_{n,i}=b_{n,n+2-i}$ for all $n\ge 1$, $2\le i\le n+1$.


\bigskip

Since we will frequently use two different recursions for $b_{n,i}$ (both mentioned in \cite{Sl}), we record them here for ease of further reference:
\begin{equation}\label{row}
b_{n+1,i}=\sum_{j=1}^{i}b_{n,j},\quad n,i\geq 1,
\end{equation}

\begin{equation}\label{twoterms}
b_{n,i+1}=b_{n,i}+b_{n-1,i+1},\quad n\geq 2, i\geq 1.
\end{equation}

Our last preliminary result concerns the enumeration of another subfamily of $321$-avoiding permutations. In the statement we will also make use of the \emph{multinomial coefficient} $\inlinemultinom{u}{v}$, which counts the number of multisets of cardinality $v$ over a set of cardinality $u$.

\begin{lemma}
\label{enumerazionefinaleprimaparte}
The number of $321$-avoiding permutations of length $n+k$ whose $k$ largest elements are LTR maxima is
\[
\sum_{i=0}^{n-1} \multinom{n-i+1}{k} b_{n,i+1}.
\]
\end{lemma}
\begin{proof}
Let $\pi'\in S_{n}(321)$. Thus, if $\pi'=M_1 P_1\cdots M_t$, we have that $P_1\cdots P_{t-1}$ is an increasing sequence. How many ways do we have to insert $n+1,\ldots,n+k$ into $\pi'$ in order to obtain a permutation $\pi\in S_{n+k}(321)$ such that the $k$ largest elements are LTR maxima? We can distinguish two cases.\\
If $\pi'$ is the identity permutation of length $n$, then we just have to insert the $k$ elements $n+1,\cdots,n+k$ in increasing order in any of the $n+1$ possible positions, and we can do that in $\inlinemultinom{n+1}{k}$ ways.\\
If $\pi'$ is not the identity permutation of length $n$, let $i$ be the position of the last element of $M_{t-1}$, with $1\le i\le n-1$. Then we can insert the $k$ elements $n+1,\ldots,n+k$ in any of the $n-i+1$ positions following $i$ in increasing order. Indeed, inserting an element in a previous position would form an occurrence of the pattern $321$. As before, this can be done in $\inlinemultinom{n-i+1}{k}$ ways. Moreover, the permutations $\pi'\in S_n(321)$ such that the last element of $M_{t-1}$ is in position $i$ are $\sum_{j=i+1}^{n} b_{j,i}$. In fact, by removing $M_t$ from $\pi'$ we obtain a permutation of length $j$, with $i<j\le n$, that has its maximum in position $i$, and the number of these permutations is $b_{j,i}$ by definition.
Summing up, the number of $321$-avoiding permutations of length $n+k$ whose $k$ largest elements are LTR maxima is 
\begin{equation}
\label{multinomialeq}
\multinom{n+1}{k} + \sum_{i=1}^{n-1} \multinom{n-i+1}{k} \sum_{m=i+1}^{n} b_{m,i}.
\end{equation}
Observe that
\begin{equation}
\label{columnrecurrence}
\sum_{m=i+1}^{n} b_{m,i}=b_{n,i+1}.
\end{equation}
Indeed, by iteration of recurrence (\ref{twoterms}), we obtain $b_{n,i+1}=b_{n,i}+b_{n-1,i+1}=b_{n,i}+b_{n-1,i}+b_{n-2,i+1}=\dots=\sum_{m=i+1}^{n} b_{m,i}$.
Plugging (\ref{columnrecurrence}) into expression (\ref{multinomialeq}), we thus obtain
\[
\multinom{n+1}{k} + \sum_{i=1}^{n-1} \multinom{n-i+1}{k} b_{n,i+1}=\sum_{i=0}^{n-1} \multinom{n-i+1}{k} b_{n,i+1},
\]
which is the thesis.
\end{proof}

We are now ready to state our main result concerning the enumeration of preimages of permutations of the form $\pi =M_1 P_1 M_2$.

\begin{theorem}
\label{MPM_theorem}
Let $\pi=M_1 P_1 M_2\in S_n$, with $M_2 \neq \emptyset$. Then
{\small
\begin{equation}\label{MPM}
\abs{q^{-1}(\pi)}=\sum_{i=1}^{m_2} \sum_{j=0}^{i-1}\binom{i-1}{j}
\left( \sum_{l=0}^{m_1-2} \multinom{m_1-l}{j+1} b_{m_1-1,l+1} \right)
\left( \sum_{k=2}^{m_2-i+1} g_{m_2-i,k} \multinom{k}{p_1+i-j-1} \right),
\end{equation}
}
where all the summations are set to be 1 whenever the set of indices is empty. 
\end{theorem}

\begin{proof}
Our goal is to describe how to obtain a generic preimage $\sigma$ from $\pi$ in such a way that we are able to count them. Let $\mu_1$ be the last element of $M_1$.
By Lemma \ref{oltrepassare}, there is at least one element of $M_2$ which is to the left of $\mu_1$ in $\sigma$. Let $\beta$ be the rightmost of such elements, and suppose that $\beta$ is the $i$-th element of $M_2$. Clearly $1\leq i\leq m_2$. Let $j$ be the number of elements of $M_2$ smaller than $\beta$ and to the left of $\mu_1$ in $\sigma$, with $0\leq j\leq i-1$. They can be chosen in $\binom{i-1}{j}$ different ways.  
In other words, $\sigma$ can be written as $\sigma=L\mu_1 R$, where $L$ is a suitable permutation of the elements of $M_1$ other that $\mu_1$ and the $j+1$ elements of $M_2$ mentioned above, and $R$ is a suitable permutations of the remaining elements (that is, the elements of $P_1$ and the remaining elements of $M_2$). We now wish to characterize (and count) the permutations $L$ and $R$ giving rise to a preimage of $\pi$.

Concerning $L$, this has to be a permutation whose $j+1$ aforementioned elements are LTR maxima and such that $q(L)$ is an increasing sequence. This means that $L$ is a 321-avoiding permutation of length $m_1 +j$ whose largest $j+1$ elements are LTR maxima. By Lemma \ref{enumerazionefinaleprimaparte}, the number of such permutations is
\[
\sum_{l=0}^{m_1-2} \multinom{m_1-l}{j+1} b_{m_1-1,l+1}.
\]

Concerning $R$, we can construct it as follows. Start by taking a permutation $\rho$ of the $m_2 -i$ elements of $M_2$ to the right of $\beta$; notice that $\rho$ has to be a 321-avoiding permutation, since applying \texttt{Queuesort} to it returns an increasing permutation. Next insert the remaining $p_1 +i-j-1$ elements into suitable positions and preserving the order they have in $\pi$. Specifically, such elements cannot be inserted to the right of the leftmost non-LTR maximum of $\rho$, whereas any other position is allowed. This is due to the fact that, applying \texttt{Queuesort} to any permutation, the relative order of the non-LTR maxima is preserved. Therefore, denoting with $k$ the position of the leftmost non-LTR maximum of $\rho$, we have $2\leq k\leq m_2 -i$, and there are $k$ allowed positions. Notice that, if $\rho$ is the identity permutation, then the above argument cannot apply; in such a case, the total number of allowed positions is $m_2 -i+1$. Summing up, and recalling the definition of $g_{n,i}$, the total number of possible permutations $R$ is:
\[
\multinom{m_2-i+1}{p_1+i-j-1}+\sum_{k=2}^{m_2-i} g_{m_2-i,k} \multinom{k}{p_1+i-j-1}=
\sum_{k=2}^{m_2-i+1} g_{m_2-i,k} \multinom{k}{p_1+i-j-1}.
\]

Combining the contributions coming from the above arguments, we get the thesis.
\end{proof}


Formula (\ref{MPM}), in its full generality, is rather involved, and it is not easy to apply it to effectively compute the number of preimages. However, the next lemma will allow us to simplify it.

\begin{lemma}\label{ballot_simplify}
For every $n\ge 1$, $1\le i\le n$,
\begin{equation}\label{bin_transf}
b_{n,i}=\sum_{h=1}^{i-1} \binom{n-h}{n-i} b_{i-1,h}.
\end{equation}
\end{lemma}
\begin{proof}
To prove (\ref{bin_transf}) we exploit a well-known combinatorial interpretation of the ballot number $b_{n,i}$ in terms of lattice paths starting from $(0,0)$, ending at $(n-1,i-1)$, remaining below the line $y=x$ and using north steps $N=(0,1)$ and east steps $E=(1,0)$. Specifically, each of such paths can be decomposed into its longest prefix that uses $i-2$ $E$ steps, followed by the remaining suffix. It is clear that the number of such prefixes ending at point $(i-2,h-1)$ is $b_{i-1,h}$, with $1\leq h\leq i-1$. Moreover, due to the specific kind of decomposition we have chosen, the remaining suffix can be any sequence of $n-i+1$ $E$ steps and $i-h$ $N$ steps starting with $E$: there are $\binom{n-h}{n-i}$ such sequences. Summing over $h$ gives the desired formula.
\end{proof}


\begin{theorem}
Let $\pi=M_1 P_1 M_2\in S_n$, with $M_2 \neq \emptyset$. Then
\begin{equation}\label{MPM_simple}
\abs{q^{-1}(\pi)}=\sum_{i=1}^{m_2} \sum_{j=0}^{i-1}\binom{i-1}{j} b_{m_1 +j+1,m_1}\cdot b_{m_2+p_1-j, m_2-i+1}.
\end{equation}
\end{theorem}
\begin{proof}
Looking at formula (\ref{MPM}), we start by observing that
\[
\sum_{l=0}^{m_1-2} \multinom{m_1-l}{j+1} b_{m_1-1,l+1} = \sum_{l=0}^{m_1-2} \binom{m_1+j-l}{j+1} b_{m_1-1,l+1} = b_{m_1 +j+1,m_1}.
\]
Indeed, in the first equality we just express the multinomial coefficient as a binomial, and the second equality comes from Lemma \ref{ballot_simplify}. 

Moreover, the summation in (\ref{MPM}) involving the $g_{n,i}$'s can be treated analogously by means of Lemma \ref{ballot_simplify}, just recalling the Remark following Proposition \ref{bi} and suitably modifying the index of summation (namely replacing $k$ with $h=m_2-i+2-k$).
\end{proof}

Yet another way to express formula (\ref{MPM}) comes from expanding the ballot numbers of the previous corollary in terms of Catalan numbers.

\begin{cor}
For $\pi=M_1 P_1 M_2\in S_n$, the quantity $\abs{q^{-1}(\pi)}$ can be expressed as a linear combination of Catalan numbers. More precisely, for any fixed $m_2=|M_2 |$, we have that $|q^{-1}(\pi )|$ is a linear combination of the Catalan numbers $C_{m_1},C_{m_1 +1},\ldots C_{m_1 +m_2 -1}$ with polynomial coefficients in $p_1$, i.e.:
\[
|q^{-1}(\pi )|=\sum_{t=0}^{m_2 -1}\omega_{m_2 ,t}(p_1 )C_{m_1 +t},
\]
where $\omega_{m_2 ,t}(p_1 )$ is a polynomial in $p_1$ of degree $m_2 -t-1$, for all $t$. \end{cor}

\begin{proof}
Indeed, if we write $b_{m_1 +j+1,m_1}=b_{m_1 +j+1,(m_1 +j+1)-(j+1)}$, by induction on the difference $j+1$ of the indices, using recurrence (\ref{twoterms}), we get
\begin{equation}\label{ballot_catalan}
b_{m_1 +j+1,m_1}=\sum_{h=1}^{\lfloor \frac{j+1}{2}\rfloor +1}(-1)^{h-1}\binom{j+2-h}{h-1}C_{m_1 +j+1-h}.
\end{equation}

Replacing the ballot number $b_{m_1 +j+1,m_1}$ with the above linear combination of Catalan numbers into (\ref{MPM_simple}), we can indeed express $|q^{-1}(\pi )|$ as a linear combination of Catalan numbers. The largest Catalan number occurring in such a linear combination is obtained when $h$ takes its minimum value 1 and $j$ takes its maximum value $m_2 -1$, which corresponds to $C_{m_1 +m_2 -1}$. Similarly, the smallest Catalan number is $C_{m_1}$, corresponding to the minimum value of the difference $j-h$, which is $-1$. Moreover, by using the closed form for $g_{n,i}$ found in Corollary \ref{gi} and the Remark following Proposition \ref{bi}, we observe that $b_{m_2+p_1-j, m_2-i+1}$ can be written as
\[
b_{m_2+p_1-j, m_2-i+1}=\binom{2m_2 +p_1 -i-j}{m_2 -i}\frac{p_1 -j+i}{2m_2 +p_1 -i-j} \ .
\]

For any fixed $m_2$, this is a polynomial of degree $m_2 -i$ in $p_1$. To conclude the proof, we now determine the degree of the polynomial multiplying $C_{m_1 +t}$ in the above mentioned linear combination, for any $t$ in the range $[0,m_2 -1]$. Looking at (\ref{ballot_catalan}), the Catalan number $C_{m_1 +t}$ (for fixed $t$) shows up for all $j\geq t$ in (\ref{MPM_simple}), and so also for all $i\geq t+1$. Notice that, for $i=t+1$, the coefficient of $C_{m_1 +t}$ has degree $m_2 -t-1$ in $p_1$ (since, in this case, $C_{m_1 +t}$ is obtained only when $j=t$ and $h=1$ in (\ref{ballot_catalan})), whereas for $i>t+1$ the resulting polynomials have lesser degree. 
\end{proof}

Effective enumerative results can be obtained for small values of the parameter $m_2$ in formula (\ref{MPM_simple}). For instance, when $m_2 =1$, we find the same result stated in Lemma \ref{catalan}. We are also able to get simple closed formulas when $m_2 =2$ and $m_2 =3$.

\begin{cor}
For $\pi=M_1 P_1 M_2\in S_n$, we get:
\begin{itemize}

\item $\abs{q^{-1}(\pi )}=C_{m_1 +1}+(p_1 +1)C_{m_1}$, when $|M_2 |=2$;

\item $\abs{q^{-1}(\pi )}=C_{m_1 +2}+(p_1 +1)C_{m_1 +1} +\frac{1}{2}(p_1 +1)(p_1 +4)C_{m_1}$, when $|M_2 |=3$.
\end{itemize}

\end{cor}


The above corollary, together with some further calculations, seem to suggest that $\omega_{m_2 ,t}(p_1 )=\omega_{m_2 +1,t+1}(p_1 )$, for all $m_2 ,t$. This could clearly simplify the computations needed to determine $\abs{q^{-1}(\pi )}$ when $m_2$ increases.

\section{Conclusions and further work}\label{conclusion}

In the spirit of previous work on \texttt{Stacksort}, we have investigated the preimages of the map associated with the algorithm \texttt{Queuesort}, obtaining a recursive description of all preimages of a given permutation and some enumerative results concerning the number of preimages. Our approach seems to suggest that, in some sense, the structure of the map associated with \texttt{Queuesort} is a little bit easier than that of the map associated with \texttt{Stacksort}, which allows us to obtain nicer results. For instance, we have been able to find a neat result concerning the possible cardinalities for the set of preimages of a given permutation; the same thing turns out to be much more troublesome for \texttt{Stacksort} \cite{D2}.

Our paper can be seen as a first step towards a better understanding of the algorithm \texttt{Queuesort}, which appears to be much less studied than its more noble relative \texttt{Stacksort}. Since the structure of \texttt{Queuesort} appears to be slightly simpler, it is conceivable that one can achieve better and more explicit results.
In this sense, there are many (classical and nonclassical) problems concerning \texttt{Stacksort} which could be fruitfully addressed also for \texttt{Queuesort}. 

For instance, it could be very interesting to investigate properties of the iterates of the map associated with \texttt{Queuesort}. For any natural number $n$, define the rooted tree whose nodes are the permutations of length $n$, having $id_n$ as its root and such that, given two distinct permutations $\tau ,\sigma \in S_n$, $\sigma$ is a son of $\tau$ whenever $q(\sigma )=\tau$. Studying properties of this tree could give some insight on the behavior of iterates of $q$. For instance, what can be said about the average depth of such a tree? Can we find enumerative results concerning some interesting statistic in the set of permutations having a fixed depth?

Another interesting issue could be the investigation of properties of the sets $Q_n^{(k)}$. For instance, how many permutations in $Q_n^{(k)}$ avoid a given pattern $\pi$?

Moreover, following \cite{CCF}, we could consider devices consisting of two queues (both with bypass) in series, where the content of the first queue is constraint to avoid some pattern. What can we say about permutations that are sortable by such devices?

Finally, Defant discovered some surprising connections between \texttt{Stacksort} and free probability theory \cite{D3}. Can we find anything similar for \texttt{Queuesort}?

\end{document}